\documentclass[a4paper]{article}

\usepackage{natbib}

\newenvironment{kvothe}{\small\begin{quote}\begin{em}}{\end{em}\end{quote}\normalsize\vglue -0.2truecm}

\begin{document}

\title{``Not only defended but also applied'':\\The perceived absurdity of Bayesian inference}

\author{
  {\sc Andrew Gelman}\\{\em Department of Statistics and Department of Political Science,}\\ {\em Columbia University}\\ 
  {\sf gelman@stat.columbia.edu}\\
  {\sc Christian P.~Robert}\\{\em Universit\'e Paris-Dauphine, CEREMADE, IUF, and CREST, Paris}\\{\sf
xian@ceremade.dauphine.fr}
}       

\maketitle

\begin{abstract}
The missionary zeal of many Bayesians of old has been matched, in the other direction, by an attitude among some
theoreticians that Bayesian methods were absurd---not merely misguided but obviously wrong in principle.   We
consider several examples, beginning with Feller's classic text on probability theory and continuing with more
recent cases such as the perceived Bayesian nature of the so-called doomsday argument. We analyze in this note
the intellectual background behind various misconceptions about Bayesian statistics, without aiming at a complete
historical coverage of the reasons for this dismissal.
\end{abstract}


\section{A view from 1950}

Younger readers of this journal may not be fully aware of the passionate battles over Bayesian inference among
statisticians in the last half of the twentieth century.  During this period, the missionary zeal of many
Bayesians was matched, in the other direction, by a view among some theoreticians that Bayesian methods are
absurd---not merely misguided but obviously wrong in principle.   Such anti-Bayesianism could hardly be
maintained in the present era, given the many recent practical successes of Bayesian methods.  But by examining
the historical background of these beliefs, we may gain some insight into the statistical debates of today.

We begin with a {\em Note on Bayes' rule} that appeared in William Feller's classic probability text:

\begin{kvothe}
``Unfortunately Bayes' rule has been somewhat discredited by metaphysical applications of the type described
above.  In routine practice, this kind of argument can be dangerous.  A quality control engineer is concerned
with one particular machine and not with an infinite population of machines from which one was chosen at
random.  He has been advised to use Bayes' rule on the grounds that it is logically acceptable and corresponds
to our way of thinking.  Plato used this type of argument to prove the existence of Atlantis, and philosophers
used it to prove the absurdity of Newton's mechanics.  In our case it overlooks the circumstance that the
engineer desires success and that he will do better by estimating and minimizing the sources of various types
of errors in predicting and guessing.  The modern method of statistical tests and estimation is less intuitive
but more realistic.  It may be not only defended but also applied.'' --- W. Feller, \citeyear{feller:1950} (pp.\ 124--125 of the 1970 edition).\end{kvothe}

Feller believed that Bayesian inference could be {\em defended} (that is, supported via theoretical
argument) but not {\em applied} to give reliable answers to problems in science or engineering, a claim that
seems quaint in the modern context of Bayesian methods being used in problems from genetics, toxicology, and
astronomy to economic forecasting and political science.  As we discuss below, what struck us about Feller's
statement was not so much his stance as his apparent certainty.

One might argue that, whatever the merits of Feller's statement today, it might have been true back in 1950.
Such a claim, however, would have to ignore, for example, the success of Bayesian methods by Turing and others
in codebreaking during the Second World War, followed up by expositions such as \cite{good:1950}, as well as
Jeffreys's {\em Theory of Probability}, which came out in 1939.  Consider this recollection from physicist and
Bayesian E. T. Jaynes:

\begin{kvothe}
``When, as a student in 1946, I decided that I ought to learn some probability
theory, it was pure chance which led me to take the book {\em Theory of Probability} by
\citeauthor{jeffreys:1939}, from the library shelf.  In reading it, I was puzzled by something which, I am
afraid, will also puzzle many who read the present book. Why was he so much on the defensive? It seemed to me
that Jeffreys' viewpoint and most of his statements were the most obvious common sense, I could not imagine
any sane person disputing them.  Why, then, did he feel it necessary to insert so many interludes of
argumentation vigorously defending his viewpoint?  Wasn't he belaboring a straw man?  This suspicion
disappeared quickly a few years later when I consulted another well-known book on probability (Feller, 1950)
and began to realize what a fantastic situation exists in this field.  {\em The whole approach of Jeffreys was
summarily rejected as metaphysical nonsense} [emphasis added], without even a description. The author assured
us that Jeffreys' methods of estimation, which seemed to me so simple and satisfactory, were completely
erroneous, and wrote in glowing terms about the success of a `modern theory,' which had abolished all these
mistakes. Naturally, I was eager to learn what was wrong with Jeffreys' methods, why such glaring errors had
escaped me, and what the new, improved methods were.  But when I tried to find the new methods for handling
estimation problems (which Jeffreys could formulate in two or three lines of the most elementary mathematics),
I found that the new book did not contain them.'' --- E. T. Jaynes (\citeyear{jaynes:1974}).
\end{kvothe}

To return to Feller's perceptions in 1950, it would be accurate, we believe, to refer to
Bayesian inference as being an undeveloped subfield in statistics at that time, with Feller being one of many
academics who were aware of some of the weaker Bayesian ideas but not of the good stuff. This goes even without
mentioning Wald's complete class results of the 1940s. (Wald's {\em Statistical Decision Functions} got
published in \citeyear{wald:1950}.)

It is in that spirit that we consider \citeauthor{feller:1970}'s notorious dismissal of Bayesian statistics,
which is exceptional not in its recommendation---after all, as of 1950 (when the first edition of his wonderful
book came out) or even 1970 (the year of his death), Bayesian methods were indeed out of the mainstream of
American statistics, both in theory and in application---but rather in its intensity.  Feller combined a
perhaps-understandable skepticism of the wilder claims of Bayesians with a na{\"\i}ve (in retrospect) faith in
the classical Neyman-Pearson theory to solve practical problems in statistics.

To say this again:  Feller's real error was not his anti-Bayesianism (an excusable position, given that many
researchers at that time were apparently unaware of modern applied Bayesian work) but rather his casual,
implicit, unthinking belief that classical methods could solve whatever statistical problems might come up.  In
short, Feller was defining Bayesian statistics by its limitations while crediting the Neyman-Pearson theory
with the 1950 equivalent of vaporware:  the unstated conviction that, having solved problems such as inference
from the Gaussian, Poisson, binomial, etc., distributions, that it would be no problem to solve all sorts of
applied problems in the future. Indeed, we take Feller's statement about ``estimating and minimizing the
sources of various types of errors'' to be a reference to the type 1 and type 2 errors of Neyman-Pearson
theory, given that he immediately follows with an allusion to ``the modern method of statistical tests and
estimation.'' In retrospect, \citeauthor{feller:1971} was wildly optimistic that the principle of ``estimating
and minimizing the sources of various types of errors'' would continue to be the best approach to solving
engineering problems. (Feller's appreciation of what a statistical problem is seems rather moderate: the two
examples Feller concedes to the Bayesian team are (i) finding the probability a family has one child given that
it has no girl and (ii) urn models for stratification/spurious contagion, problems that are purely
probabilistic, no statistics being involved.)  Or, to put it another way, even within the context of prediction
and minimizing errors, why be so sure that Bayesian methods cannot apply?  Feller perhaps leapt from the
existence of philosophical justification of Bayesian inference, to an assumption that philosophical arguments
were the {\em only} justification of Bayesian methods.

Where was this coming from, historically?  With Stephen Stigler out of the room, we are reduced to speculation
(or, maybe we should say, we are free to speculate).  We doubt that Feller came to his own considered judgment
about the relevance of Bayesian inference to the goals of quality control engineers.  Rather, we suspect that
it was from discussions with one or more statistician colleagues that he drew his strong opinions about the
relative merits of different statistical philosophies.  In that sense, Feller is an interesting case in that he
was a leading mathematician of his area, a person who one might have expected would be well informed about
statistics, and the quotation reveals the unexamined assumptions of his colleagues.  It is doubtful that even
the most rabid anti-Bayesian of 2010 would claim that Bayesian inference cannot apply.  (We would further
argue that the ``modern methods of statistics''  Feller refers to have to be understood in historical context
as eliminating older approaches by Bayes, Laplace, and other 19th century authors, in a spirit akin to Keynes
(1921). Modernity starts with the great anti-Bayesian Ronald Fisher who, along with Richard von Mises, is
mentioned on page 6 by Feller as the originator of ``the statistical attitude towards probability.'' 
\cite{vonmises:1957} may have been strong in mathematics and other fields, but when it came to a simple
comparison of binomial variances, he didn't know how to check for statistical significance; see
\cite{gelman:2010}. He rejected not only ``persistent subjectivists'' (p.\ 94) such as John Maynard Keynes and
Harold Jeffreys, but also Fisher's likelihood theory (p.\ 158).)

\section{The link between Bayes and bogosity}

Non-Bayesians still occasionally dredge up Feller's quotation as a pithy reminder of the perils of philosophy
unchained by empiricism (see, for example, \citealp{ryder:1976}, and \citealp{diNardo:2008}).  In a recent
probability text, \cite{stirzaker:1999} reviews some familiar probability paradoxes (e.g., the Monty Hall
problem) and draws the following lesson:

\begin{kvothe}
``In any experiment, the procedures and rules that define the sample space and all the probabilities must be explicit and fixed before you begin.  This predetermined structure is called a protocol.  Embarking on experiments without a complete protocol has proved to be an extremely convenient method of faking results over the years.  And will no doubt continue to be so.''
\end{kvothe}

Strirzaker follows up with a portion of the Feller quote and writes, ``despite all this experience, the popular
press and even, sometimes, learned journals continue to print a variety of these bogus arguments in one form or
another.''  We are not quite sure why he attributes these problems to Bayes, rather than, say, to
Kolmogorov---after all, these error-ridden arguments can be viewed as misapplications of probability theory
that might never have been made if people were to work with absolute frequencies rather than fractional
probabilities \citep{vonmises:1957,gigerenzer:2002}.

In any case, no serious scientist can be interested in bogus arguments (except, perhaps, as a teaching tool or
as a way to understand how intelligent and well-informed people can make evident mistakes, as discussed in
chapter 3 of  \cite{gelman:park:etal:2008}).  What is perhaps more interesting is the presumed association
between Bayes and bogosity.  We suspect that it is Bayesians' openness to making assumptions that makes their
work a particular target, along with  (some) Bayesians' intemperate rhetoric about optimality.  Somehow
classical terms such as ``uniformly most powerful test'' do not seem so upsetting.
%, perhaps because it is clear that said uniformity relies on assumptions which are so clearly invalid (in most
%practical circumstances) as to be considered idealizations at best. However, this was not the case for Feller,
%who genuinely seemed to think that the optimality properties of ``modern methods of statistical test and
%estimation" qualified them as optimal in practice.
Perhaps what has bothered mathematicians such as Feller and Stirzaker is that Bayesians actually seem to
believe their assumptions rather than merely treating them as counters in a mathematical game. In the first
quote, the interpretation of the prior distribution as a reasoning based on an ``infinite population of
machines'' certainly indicates that Feller takes the prior at face value! As shown by the recent foray of
\cite{burdzy:2009} into the philosophy of Bayesian foundations and in particular of deFinetti's, this
interpretation may be common among probabilists, whereas we see applied statisticians as considering both prior
and data models as assumptions to be valued for their use in the construction of effective statistical
inferences.

In applied Bayesian inference, it is not necessary for us to believe our assumptions, any more than
biostatisticians believe in the truth of their logistic regressions and proportional hazards models.  Rather,
we make strong assumptions and use subjective knowledge in order to make inferences and predictions that can be
tested by comparing to observed and new data (see \citealp{gelman:shalizi:2010}, or \citealp{mayo:1996} for a
similar attitude coming from a non-Bayesian direction). Unfortunately, we doubt Stirzaker was aware of this
perspective when writing his book---nor was Feller,  working years before either of the present authors were
born.

Recall the following principle, to which we (admitted Bayesians) subscribe:
\begin{kvothe}
Everyone uses Bayesian inference when it is clearly appropriate.  A Bayesian is someone who uses Bayesian inference even when it might seem inappropriate.
\end{kvothe}
What does this mean?  Mathematical modelers from R. A. Fisher on down have used and will use probability to
model physical or algorithmic processes that seem well-approximated by randomness, from rolling of dice to
scattering of atomic particles to mixing of genes in a cell to random-digit dialing.  To be honest, most
statisticians are pretty comfortable with probability models even for processes that are not so clearly
probabilistic, for example fitting logistic regressions to purchasing decisions or survey responses or
connections in a social network.  (As discussed in \citealp{robert:2010}, Keynes' {\em Treatise on Probability}
is an exception in that \citeauthor{keynes:1921} even questions the sampling models.) Bayesians will go the
next step and assign a probability distribution to a parameter that one could not possibly imagine to have been
generated by a random process, parameters such as the coefficient of party identification in a regression on
vote choice, or the overdispersion in a network model, or Hubble's constant in cosmology. There is no
inconsistency in this opposition once one realizes that priors are not reflections of a hidden ``truth'' but
rather evaluations of the modeler's uncertainty about the parameter. Using distributions on a fixed but
unknown parameter extends to non-Bayesians like \cite{efron:1986} and \cite{fraser:2011}.

As noted above, it is our impression that the assumptions of the likelihood are generally more crucial---and
often less carefully examined---than the assumptions in the prior.  Still, we recognize that Bayesians take
this extra step of mathematical modeling.  In some ways, the role of Bayesians compared to other statisticians
is similar to the position of economists compared to other social scientists, in both cases making additional
assumptions that are clearly wrong (in the economists' case, models of rational behavior) in order to get
stronger predictions.  With great power comes great responsibility, and Bayesians and economists alike have the
corresponding duty to check their predictions and abandon or extend their models as necessary.

To return briefly to Stirzaker's quote, we believe he is wrong---or, at least, does not give any good
evidence---in his claim that ``in any experiment, the procedures and rules that define the sample space and all
the probabilities must be explicit and fixed before you begin.'' Setting a protocol is fine if it is practical,
but as discussed by Rubin (1976), what is really important from a statistical perspective is that all the
information used in the procedure be based on known and measured variables.  This is similar to the idea in
survey sampling that clean inference can be obtained from probability sampling---that is, rules under which all
items have nonzero probabilities of being selected, with these probabilities being known (or, realistically,
modeled in a reasonable way).

It is unfortunate that certain Bayesians have published misleading and oversimplified expositions of the Monty
Hall problem (even when fully explicated, the puzzle is not trivial, as the resolution requires a full
specification of a probability distribution for Monty's possible actions under various states of nature, see
e.g.~\citealp{rosenthal:2006}); nonetheless, this should not be a reason for statisticians to abandon decades
of successful theory and practice on adaptive designs of experiments and surveys, not to mention the use of
probability models for non-experimental data (for which there is no ``protocol'' at all).

\section{The sun'll come out tomorrow}

The prequel to Feller's quotation above is the notorious argument, attributed to Laplace, that uses a flat
prior distribution on a binomial probability to estimate the probability the sun will rise tomorrow.  The idea
is that the sun has risen $n$ out of $n$ successive days in the past, implying a posterior mean of
$(n+1)/(n+2)$ of the probability $p$ of the sun rising on any future day. (\citealp{gorroochurn:2011} gives a
recent coverage of the many criticisms that ridiculed Laplace's ``mistake.'')

To his credit, Feller immediately recognized the silliness of that argument.  For one thing, we don't have
direct information on the sun having risen on any particular day, thousands of years ago, and cannot predict
what will occur the next morning.  So the analysis is conditioning on data that don't exist, in the sense that
the assumed model is not supported by the actual evidence.

More than that, though, the big, big problem with the Pr(sunrise tomorrow $|$ sunrise in the past) argument is
not in the prior but in the likelihood, which assumes a constant probability and independent events.  Why
should anyone believe that?  Why does it make sense to model a series of astronomical events as though they
were spins of a roulette wheel in Vegas?  Why does stationarity apply to this series? That's not frequentist,
it isn't Bayesian, it's just dumb.  Or, to put it more charitably, it's a plain vanilla default model that we
should use only if we are ready to abandon it on the slightest pretext. The Laplace law of
succession has been discussed {\em ad nauseaum} in relation to the Humean debate about inference (see, e.g.,
\citealp{sober:2008}).  Furthermore, \cite{berger:bernardo:sun:2009} discuss other prior distributions for the
model.  Here, however, we are focusing on the likelihood function, which, despite its extreme inappropriateness
for this problem, is typically accepted without question.

It is no surprise that when this model fails, it is the likelihood rather than the prior that is causing the
problem.  In the binomial model under consideration here, the prior comes into the posterior distribution only
once, and the likelihood comes in $n$ times.  It is perhaps merely an accident of history that skeptics and
subjectivists alike strain on the gnat of the prior distribution while swallowing the camel that is the
likelihood.  In any case, it is instructive that Feller saw this example as an indictment of Bayes (or at least
of the uniform prior as a prior for ``no advance knowledge'') rather than of the binomial distribution.

\section{The ``doomsday argument'' and confusion between frequentist and Bayesian ideas}

Bayesian inference has such a hegemonic position in philosophical discussions that, at this point, statistical
arguments get interpreted as Bayesian even when they are not.  

An example is the so-called doomsday argument \citep{carter:mccrea:1983}, which holds that there is a high probability that
humanity will be extinct (or will drastically reduce in population) soon, because if this were not true---if, for
example, humanity were to continue with 10 billion people or so for the next few thousand years---then each of
us would be among the first people to exist, and that's highly unlikely.  To put it slightly more formally, the
``data'' here is the number of people, $n$, who have lived on Earth up to this point, and the ``hypothesis''
corresponds to the total number of people, $N$, who will ever live.  The statistical argument is that $N$ is
almost certainly within two orders of magnitude of $n$, otherwise the observed $n$ would be highly improbable.
And if $N$ cannot be much more than $n$, this implies that civilization cannot exist in its current form for
millenia to come.

For our purposes here, the (sociologically) interesting thing about this argument is that it's been presented
as Bayesian (see, for example, \citealp{dieks:1992}) but it isn't a Bayesian analysis at all!  The ``doomsday
argument'' is actually a classical frequentist confidence interval.  Averaging over all members of the group
under consideration, 95\% of these confidence intervals will contain the true value.  Thus, if we go back and
apply the doomsday argument to thousands of past data sets, its 95\% intervals should indeed have 95\%
coverage.  In 95\% of populations examined at a randomly-observed rank, $n$ will be between $0.025N$ and
$0.975N$.  This is the essence of Neyman-Pearson theory, that it makes claims about averages, not about
particular cases.

However, this does not mean that there is a 95\% chance that any particular interval will contain the true
value.  Especially not in this situation, where we have additional subject-matter knowledge.  That's where
Bayesian statistics (or, short of that, some humility about applying classical confidence statements to
particular cases) comes in.  The doomsday argument seems silly to us, and we see it as fundamentally not
Bayesian. Some Bayesian versions of the doomsday argument have been constructed, but, from our
perspective, these are just unsuccessful attempts to take what is fundamentally a frequentist idea and adapt it
to make statements about particular cases.  See \cite{dieks:1992} and \cite{neal:2006} for detailed critiques
of the assumptions underlying Bayesian formulations of the doomsday argument.

The doomsday argument sounds Bayesian, though, having three familiar features 
that are (unfortunately) sometimes associated with traditional Bayesian reasoning:
\begin{itemize}
\item It sounds more like philosophy than science.
\item It's a probabilistic statement about a particular future event.
\item It's wacky, in an overconfident, ``you gotta believe this counterintuitive finding, 
it's supported by airtight logical reasoning,'' sort of way.  
\end{itemize}
Really, though, it's a classical confidence interval, tricked up with enough philosophical mystery and
invocation of Bayes that people think that the 95\% interval applies to every individual case.  Or, to put it
another way, the doomsday argument is the ultimate triumph of the idea, beloved among Bayesian educators, that
our students and clients don't really understand Neyman-Pearson confidence intervals and inevitably give them
the intuitive Bayesian interpretation.

Misunderstandings of the unconditional nature of frequentist probability statements are hardly new.  Consider
Feller's statement, ``A quality control engineer is concerned with one particular machine and not with an
infinite population of machines from which one was chosen at random.''  It sounds as if Feller is objecting to
the prior distribution or ``infinite population,'' $p(\theta)$, and saying that he only wants inference for a
particular value of $\theta$. This misunderstanding is rather surprising when issued by a probabilist but it
shows a confusion between data and parameter: as mentioned above, the engineer wants to condition upon the data
at hand (with a specific if unknown value of $\theta$ lurking in the background). Again,
this relates to Feller holding a second-hand opinion on the topic and backing it with a cooked-up story.  It
does not help that many Bayesians over the years have muddied the waters by describing parameters as random
rather than fixed.  Once again, for Bayesians as much as for any other statistician, parameters are (typically) fixed but
unknown.  It is the knowledge about these unknowns that Bayesians model as random.

In any case, we suspect that many quality control engineers do take measurements on multiple machines, maybe even
populations of machines, but to us Feller's sentence noted above has the interesting feature that it is
actually the opposite of the usual demarcation:  typically it is the Bayesian who makes the claim for inference
in a particular instance and the frequentist who restricts claims to infinite populations of replications.

\section{Conclusions}

Why write an article picking on sixty years of confusion?  We are not seeking to malign the reputation of
Feller, a brilliant mathematician and author of arguably the most innovative and intellectually stimulating
book ever written on probability theory.  Rather, it is Feller's
brilliance and eminence that makes his attitude that much more interesting:  that this centrally-located figure in probability theory could make a
statement that could seem so silly in retrospect (and even not so long in retrospect, as indicated by the
memoir of Jaynes quoted above).

Misunderstandings of Bayesian statistics can have practical consequences in the present era as well.  We could
well imagine a reader of Stirzaker's generally excellent probability text taking home the message that
all probabilities ``must be explicit and fixed before you begin,'' thus missing out on some of the most exciting
and important work being done in statistics today.

In the last half of the twentieth century, Bayesians had the reputation (perhaps deserved) of being
philosophers who were all too willing to make broad claims about rationality, with optimality theorems that were
ultimately built upon questionable assumptions of subjective probability, in a denial of the
garbage-in-garbage-out principle, thus defying all common sense.  In opposition to this nonsense, Feller (and
others of his time) favored a mixture of Fisher's rugged empiricism and the rigorous Neyman-Pearson theory, which ``may be not only defended but also applied.''  And, indeed, if the classical theory of hypothesis testing had lived up to the promise it seemed to have in 1950 (fresh after solving important operations-research problems in the Second World War), then indeed maybe we could have stopped right there.

But, as the recent history of statistics makes so clear, no single paradigm---Bayesian or otherwise---comes
close to solving all our statistical problems (see the recent reflections of \cite{senn:2011}) and there are huge
limitations to the type-1, type-2 error framework which seemed so definitive to Feller's colleagues at the
time.  At the very least, we hope Feller's example will make us wary of relying on the advice of colleagues to
criticize ideas we do not fully understand.  New ideas by their nature are often expressed awkwardly and with
mistakes---but finding such mistakes can be an occasion for modifying and improving these ideas rather than
rejecting them.

\section*{Acknowledgements} 
We thank David Aldous, Ronald Christensen, the Associate Editor, and two reviewers for helpful comments.  In
addition, the first author (AG) thanks the Institute of Education Sciences, Department of Energy, National
Science Foundation, and National Security Agency for partial support of this work.  He remembers reading with
pleasure much of Feller's first volume in college, after taking probability but before taking any statistics
courses.  The second author's (CPR) research is partly supported by the Agence Nationale de la Recherche (ANR,
212, rue de Bercy 75012 Paris) through the 2007--2010 grant ANR-07-BLAN-0237 ``SPBayes.'' He remembers buying
Feller's first volume in a bookstore in Ann Arbor during a Bayesian econometrics conference where he was kindly
supported by Jim Berger.
%\end{acknowledgments}

\renewcommand{\bibsection}{\section*{References}}

%\bibliographystyle{ba}
%\bibliography{biblio.bib}

\end{document}